\documentclass[12pt]{article}
\usepackage{amssymb,amsmath, amsthm}
\usepackage{graphicx}

\newtheorem{thm}{Theorem}

\newtheorem{cor}{Corollary}

\newtheorem{lemma}{Lemma}

\newenvironment{defin}{\medskip\noindent{\sc
Definition}. }{\goodbreak\medskip}
\newenvironment{nota}{\medskip\noindent{\sc
Notations}. }{\goodbreak\medskip}
\newenvironment{remk}{\noindent{\sc
Remarks}. }{\goodbreak\vskip10pt}
\newtheorem{prop}{Proposition}

\newtheorem{ques}{Question}
\newtheorem{exa}{Example}

\def\demo{\medskip\goodbreak\noindent
     \hbox{\sc Proof \kern .3em}\ignorespaces}%
  \def \qedbox{$\square$}%
  \def \qed{\hglue1mm\hfill{\ifmmode\qedbox
     \else\unskip\ \hglue0mm\hfill\qedbox\medskip
      \goodbreak\fi}}%
\def\enddemo{\qed\goodbreak\vskip10pt}%

\def\qed{\hglue1mm\hfill\raise -2pt\hbox{\vrule\vbox to 10pt{\hrule width
4pt
                  \vfill\hrule}\vrule}}

\newcommand{\U}{\mathbb {U}}
\newcommand{\T}{\mathbb {T}}

\newcommand{\A}{\mathbb {A}}
\newcommand{\esse}{\mathbb {S}}
\newcommand{\R}{\mathbb {R}}

\newcommand{\C}{\mathbb {C}}
\newcommand{\Z}{\mathbb {Z}}
\newcommand{\N}{\mathbb {N}}
\newcommand{\Nc}{\mathcal {N}}

\newcommand{\Dc}{\mathcal {D}}
\newcommand{\Hc}{\mathcal {H}}

\newcommand{\Gc}{\mathcal {G}}
\newcommand{\Lc}{\mathcal {L}}
\newcommand{\Fc}{\mathcal {F}}

\newcommand{\Sc}{\mathcal {S}}

\begin{document}
\title{Lyapunov exponents for conservative twisting dynamics: a survey}
\author{M.-C. ARNAUD 
\thanks{ANR-12-BLAN-WKBHJ}
\thanks{Avignon Universit\'e , Laboratoire de Math\'ematiques d'Avignon (EA 2151),  F-84 018 Avignon,
France. e-mail: Marie-Claude.Arnaud@univ-avignon.fr} 
\thanks{membre de l'Institut universitaire de France}}
\maketitle
\abstract{\noindent \sl  Finding special orbits (as periodic orbits) of dynamical systems by variational methods and especially by minimization methods is an old method (just think of the geodesic flow). More recently,  new results concerning the existence of minimizing sets and minimizing measures were proved in the setting of conservative twisting dynamics.  These twisting dynamics include geodesic flows as well as the dynamics close to a completely elliptic periodic point of a symplectic diffeomorphism where the torsion is positive definite (this implies the existence of a normal form $(\theta, r)\mapsto (\theta+\beta r+o(r), r+o(r))$ with $\beta$ positive definite). Two aspects of this theory are called the Aubry-Mather theory and the weak KAM theory. They were built by Aubry \& Mather in the '80s in the 2-dimensional case and by Mather, Ma\~n\'e and Fathi in the '90s in higher dimension. 

We will explain what are the conservative  twisting dynamics and summarize the existence results of minimizing measures. Then we will explain more recent results concerning the link between different notions for  minimizing measures for twisting dynamics:
\begin{enumerate}
\item[$\bullet$] their Lyapunov exponents;
\item[$\bullet$] their Oseledets splitting;
\item[$\bullet$] the shape of the support of the measure.
\end{enumerate}
The main question in which we are interested is:  given some minimizing measure of a conservative twisting dynamics, is there a link between the geometric shape of its support and its Lyapunov exponents? Or : can we deduce the Lyapunov exponents of the measure  from the ``shape'' of the support of this measure? \\
Some proofs but not all of them will be provided. Some questions are raised in the last section.

}

\newpage

\tableofcontents

  \newpage

   \section{Twisting conservative dynamics}\label{sectiontwist}
   All the dynamics we study here are defined on the cotangent bundle $T^*M$ of some closed manifold $M$, endowed with its usual symplectic form $\omega$. More precisely, if $q=(q_1, \dots, q_n) $ are some coordinates on $M$, we complete them with  their  dual coordinates $p=(p_1, \dots , p_n)$ to obtain some coordinates on $T^*M$: if $\lambda\in T^*M$ is a  $1$-form on $M$, then  its coordinates $p_1,\dots , p_n$ are given by  $\displaystyle{\lambda=\sum_{i=1}^n p_idq_i}$. The expression of the  symplectic form in these coordinates is $\displaystyle{\omega=dq\wedge dp=\sum_{i=1}^ndq_i\wedge dp_i}$.  A change of coordinates of $M$ doesn't change the symplectic form $\omega$ and then the definition is correct. We will generally use the notation $(q, p)$ for such coordinates.\\
   When $M=\T^n$, we will identify $T^*M$ with the $2n$-dimensional annulus $\A_n=\T^n\times \R^n$. \\
   Let us recall that a diffeomorphism $f$ of $T^*M$ is {\em symplectic} if it preserves the symplectic form: $f^*\omega=\omega$. 
   \subsection{A local notion: the twist condition}
   \begin{nota}
   We denote by $\pi:T^*M\rightarrow M$ the canonical projection $(q, p)\mapsto q$.\\
   At every $x=(q,p)\in T^*M$, we define the vertical subspace $V(x)=\ker D\pi(x)\subset T_x(T^*M)$ as being the tangent subspace at $x$ to the fiber $T^*_qM$.
   \end{nota}
   \begin{exa} A symplectic 
   $C^1$ diffeomorphism $f:\A_1\rightarrow \A_1$ of the 2-dimensional annulus is a {\em positive symplectic twist map} if
   \begin{enumerate}
   \item[$\bullet$] it is homotopic to identity;
   \item[$\bullet$] it twists the vertical to the positive side: $\forall x\in \A_1, D(\pi\circ f)(x).\begin{pmatrix} 0\\ 1\end{pmatrix}>0$.
   \begin{center}
\includegraphics[width=6cm]{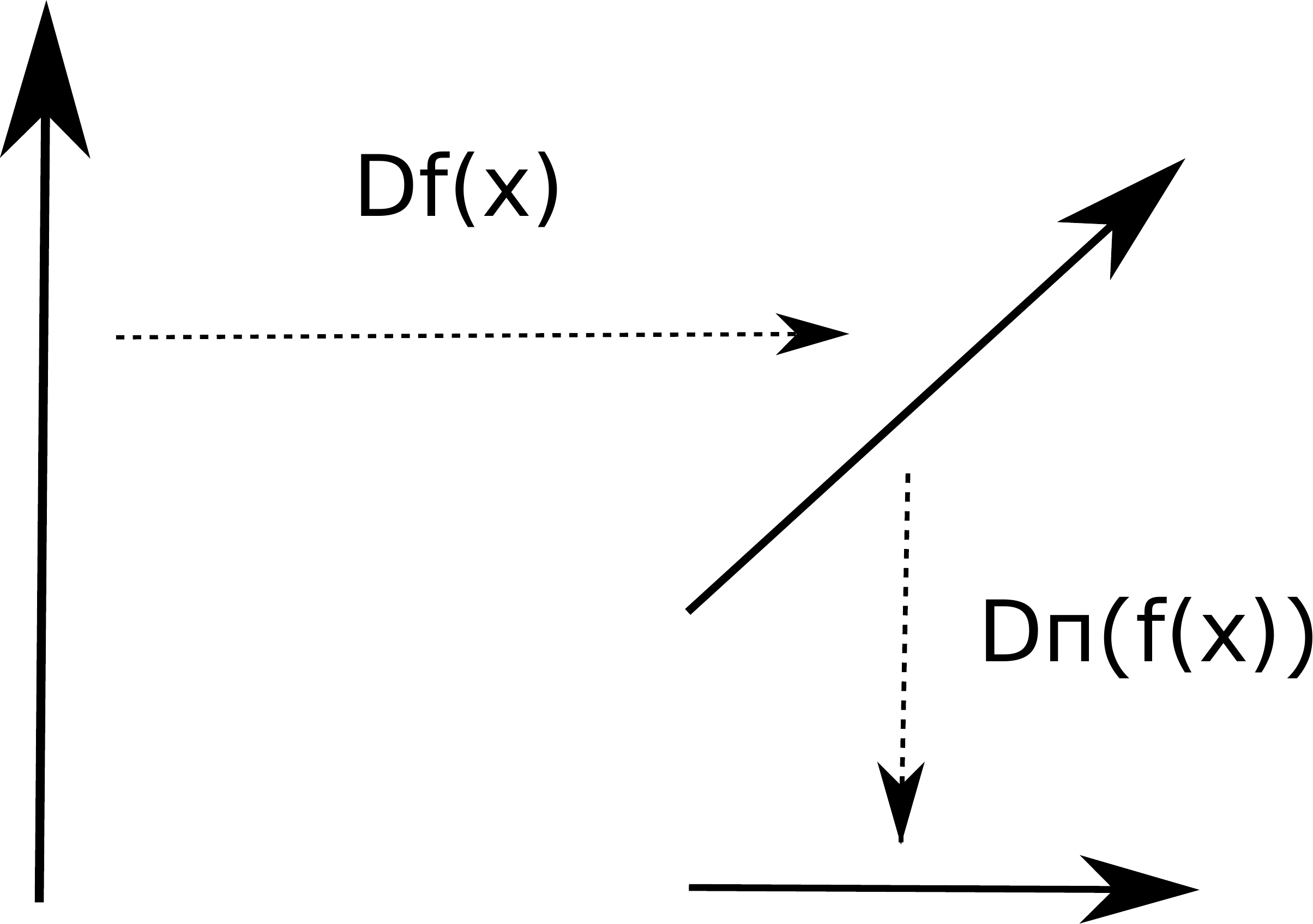}
\end{center}
  \noindent There exists of course a notion of {\em negative symplectic twist map}. 
   \end{enumerate}
   The notion of twist map that we introduce  is local, but the (local) twist condition implies a global property: when we unfold the cylinder (i.e. we are in the universal covering $\R^2$ of $\A_1$ and we consider a lift of the twist map), the image of a fiber $\{ q\}\times \R$ by the lift of a symplectic twist map is then a graph above a part of of $\R\times\{ 0\}$:
    \begin{center}
\includegraphics[width=6cm]{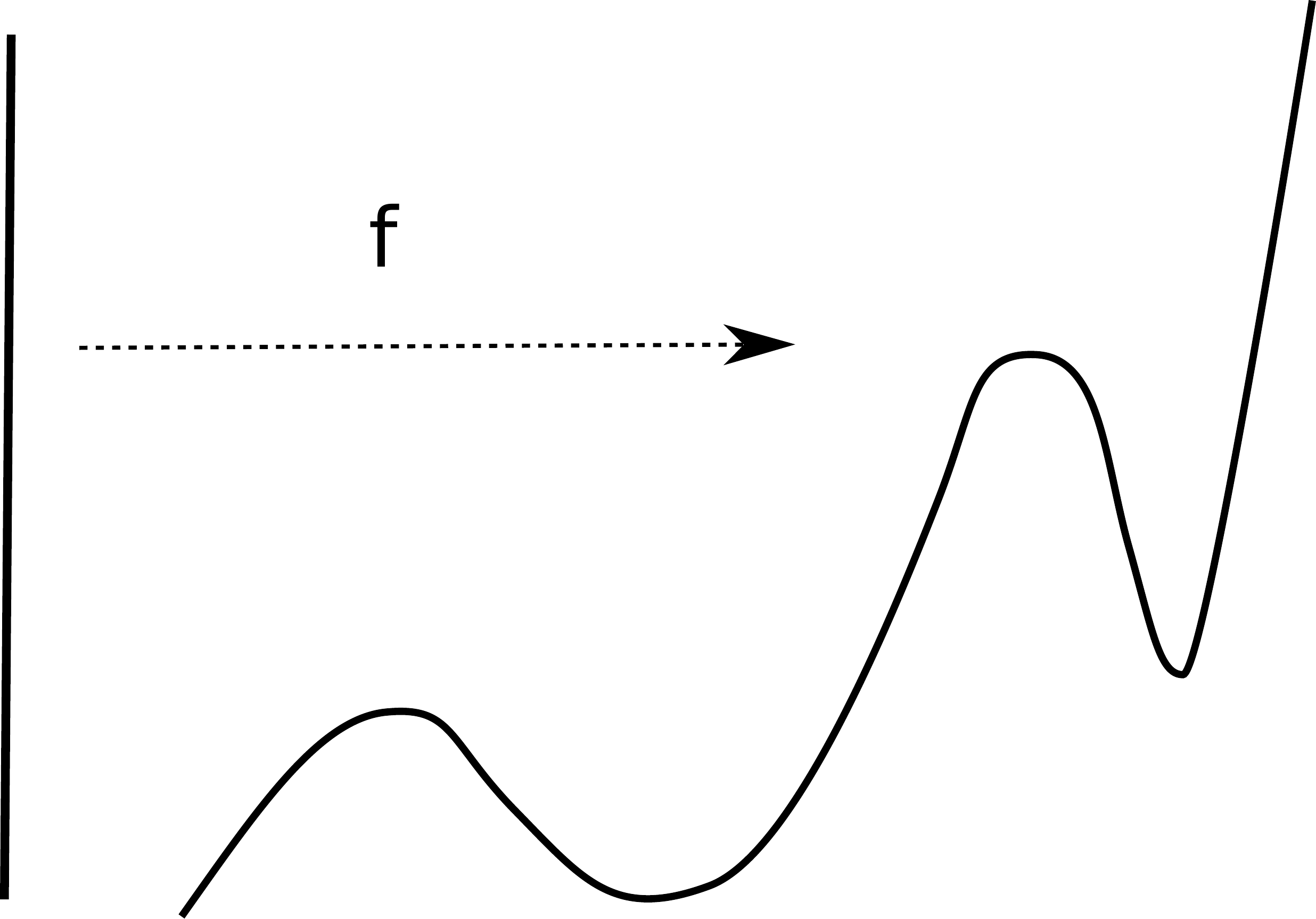}
\end{center} 
For example, the map $f_0:(q, p)\mapsto (q+p,p)$ is a symplectic twist map of $\A_1$.  \end{exa}

  When thinking of a possible extension of the notion of twisting dynamics to higher dimension,   the first possibility is to ask that the image of any vertical $V(x)$ by the tangent dynamics $Df $ is transverse to the vertical $V(f(x))$.  If we express $Df$  in two charts of coordinates $(q,p)$: 
  $$Df=\begin{pmatrix}a&b\\ c&d
  \end{pmatrix}$$
  this is equivalent to ask that $\forall x, \det b(x)  \not=0$. When $M=\T^n$, such  diffeomorphims were introduced and studied by M.~Herman in \cite{Her1}, where he called them {\em monotone}. When $\dim M\geq 2$, the monotonicity condition  doesn't imply that the image of a fiber is a graph above the zero section, even if $M=\T^n$ and if we unfold the $2n$-dimensional annulus. For example, the map $f:\A_2\rightarrow \A_2$ defined for $(q, p)\in \T^2\times \R^2=\C/(\Z+i\Z)\times \R^2$ by $f(q,p)=(q+e^{p_1-ip_2}, p)$ is monotone but the projection of the restriction of its lift  to any fiber is not injective. 
    
  If $f$ is a twist map of the 2-dimensional annulus, $f^2$ is not necessarily a twist map: the twist condition is just valuable for ``small times'' (here time 1). 
  
  Hence for Hamiltonians, we will translate the twist condition for small times. In coordinates $(q, p)$, the Hamilton equations for $H\in C^2(T^*M, \R)$ are: 
  $$\dot q=\frac{\partial H}{\partial p}(q,p);\quad \dot p=-\frac{\partial H}{\partial q}(q, p).$$
Let us denote the Hamiltonian flow of $H$ by $(\varphi_t^H)$ and let $(\delta q, \delta p)$ be an infinitesimal solution, i.e. $\begin{pmatrix} \delta q(t)\\ \delta p(t)\end{pmatrix}=D\varphi_t(q(0), p(0)).\begin{pmatrix} \delta q(0)\\ \delta p(0) \end{pmatrix}$.
By differentiating the Hamilton equations, we obtain $\delta \dot q =\frac{\partial^2 H}{\partial q \partial p}\delta q+\frac{\partial ^2 H}{\partial p^2}\delta p$ and then 
$$D(\pi\circ \varphi_t)(q(0), p(0)).\begin{pmatrix}0\\ \delta p\end{pmatrix}=t\frac{\partial ^2H}{\partial p^2}(q(0), p(0))\delta p+o(t).$$
We will say that the Hamiltonian $H$ satisfies the twist condition if at every point $\frac{\partial^2 H}{\partial p^2}$ is non-degenerate. In this case, even for small times, the Hamiltonian flow is not necessarily a twist map; indeed, the $o(t)$ above is not uniform in $(q, p)$.

  \subsection{Global notions: globally positive diffeomorphisms and Tonelli Hamiltonians}

 Unfortunately, we are able to do nothing with the local definition of twisting dynamics that we gave in the above subsection.  
 
 There are two problems:
 \begin{enumerate}
 \item \label{pt1} we need to find some special invariant subsets for the dynamics; 
  \item \label{pt2}  we want to say something about the Lyapunov exponents along these invariant subsets.
 \end{enumerate}
 In general there are two main ways to find invariant subsets for those dynamics: perturbative methods and variational methods. Perturbative methods, as K.A.M. theorems are, are valuable close to completely integrable dynamics (see \cite{Her1} for the definition in the case of the $2n$-dimensional annulus). C.~Gol\'e gives in \cite{Gol1}, section 27.B   a similar condition, that he calls ``asymptotic linearity'', that makes possible the use of variational methods in this perturbative case. But we won't explain the perturbative case in this survey. We will only work in the so-called {\em coercive} case (see section 27.B of \cite{Gol1} for example). 
 
 More precisely, we will make some assumptions such that
 \begin{enumerate}
 \item[$\bullet$] we can associate a function $\Fc$ to the dynamics, such that the critical points of $\Fc$ are in some sense the orbits for the dynamics;
 \item[$\bullet$] the function $\Fc$ admits some minima, and then  some ``minimizing orbits''. For the symplectic twist maps of the 2-dimensional annulus, these minimizing orbits are the heart of the theory that S.~Aubry and J.~Mather  independently developped at the beginning of  the '80s (see \cite{ALD} and \cite{Mat1}).
 \end{enumerate}
 
 That is why we introduce the following definitions. The first one comes from \cite{MMS} and \cite{Gol1}, the second one is very classical.  \subsubsection{Globally positive diffeomorphisms}\label{ssgpd}
 \begin{defin}
 A {\em globally positive} diffeomorphism of $\A_n$ is a symplectic $C^1$ diffeomorphism $f:\A_n\rightarrow \A_n$ that is homotopic to ${\rm Id}_{\A_n}$ and that has a lift  $F: \R^n\times \R^n\rightarrow \R^n\times \R^n$  that admits a $C^2$ generating function $S: \R^n\times \R^n\rightarrow \R$ such that:
 \begin{enumerate}
 \item[$\bullet$] $\forall k\in\Z^n, S(q+k, Q+k)=S(q, Q)$;
 \item[$\bullet$] there exists $\alpha>0$ such that: $\frac{\partial^2 S}{\partial q\partial Q}(q, Q)(v, v)\leq -\alpha\| v\|^2$;
 \item[$\bullet$] $F$ is implicitly given by:
 $$F(q,p)=(Q, P)\Longleftrightarrow \left\{\begin{matrix}p=-\frac{\partial S}{\partial q}(q, Q)\\ P=\frac{\partial S}{\partial Q}(q, Q)
 \end{matrix}\right.$$
 \end{enumerate}
  \end{defin}
  \begin{exa}\label{exa2}
  The diffeomorphism $F_0: (q,p)\in\R^n\times\R^n\mapsto (q+p, p)\in\R^n\times\R^n$ is the lift of a globally positive diffeomorphism $f_0$ of $\A_n$ and a generating function associated to $F_0$ is defined by $S_0(q, Q)=\frac{1}{2}\| q-Q\|^2$.
  \end{exa}
 If $f$,   $F$ satisfy the above hypotheses, the restriction to any fiber $\{q\}\times \R^n$ of $\pi\circ F$ and $\pi\circ F^{-1}$ are diffeomorphisms (a proof is given in \cite{Gol1}). In particular, this implies that $f$ is (locally) monotone.
 
 Moreover, for every $k\geq 2$, $q_0, q_k\in\R^n$, the function $\Fc:(\R^n)^{k-1}\rightarrow \R$ defined by $\displaystyle{\Fc(q_1, \dots, q_{k-1})=\sum_{j=1}^kS(q_{j-1}, q_j)}$ has a minimum, and at every critical point for $\Fc$, the following sequence is a piece of orbit for $F$:
 $$(q_0, -\frac{\partial S}{\partial q}(q_0, q_1)), (q_1,  \frac{\partial S}{\partial Q}(q_0, q_1)), (q_2,    \frac{\partial S}{\partial Q}(q_1, q_2)), \dots , (q_k,  \frac{\partial S}{\partial Q}(q_{k-1}, q_k)).$$
 For the map $F_0$ defined in Example \ref{exa2}, the function $\Fc_0$ defined by \\
 $\displaystyle{\Fc_0(q_1, \dots, q_{k-1})=\frac{1}{2}\sum_{j=1}^k\| q_{j-1}-q_j\|^2}$ attains its minimum at its unique critical point $(q_1, \dots, q_{k-1})=(q_0+\frac{q_k-q_0}{k}, q_0+2\frac{q_k-q_0}{k}, \dots, q_0+(k-1)\frac{q_k-q_0}{k})$ and the corresponding piece of orbit is:
 $$(q_0, \frac{q_k-q_0}{k}), (q_1, \frac{q_k-q_0}{k}), \dots , (q_k, \frac{q_k-q_0}{k}).$$

 Let us discuss a little the condition on $\frac{\partial^2 S}{\partial q\partial Q}$.  If the matrix of $Df$ in coordinates $(q,p)$ is $Df=\begin{pmatrix}a&b\\ c&d
  \end{pmatrix}$, we have $$\left( b(q , p)\right)^{-1}=-\frac{\partial^2 S}{\partial q\partial Q}(q ,Q).$$
  Hence the condition that we gave for the partial derivatives of $S$ can be rewritten in terms of matrices: $b^{-1}+{}^tb^{-1}\geq \alpha{\bf 1}$ where ${\bf 1}$ is the identity matrix and we use the usual order for the symmetric matrices.

The reader could think of some other possible notions of global twist, for which ${}^tb^{-1}+b^{-1}$ is indefinite. But in this case, very pathological phenomena can occur;  M.~Herman showed very  strange phenomena   in the case of  a ``normal indefinite torsion'' in \cite{Her2} (the torsion is ${}^tb+b$ and it has the same signature as ${}^tb^{-1}+b^{-1}={}^tb^{-1}({}^tb+b)b^{-1}$).
 \subsubsection{Tonelli Hamiltonians} 
 A $C^2$ function $H: \T^*M\rightarrow \R$ is a {\em Tonelli Hamiltonian} if it is:
 \begin{enumerate}
 \item[$\bullet$] superlinear in the fiber, i.e. $\forall A\in \R, \exists B\in \R, \forall (q,p)\in T^*M, \| p\|\geq B\Rightarrow  H(q,p) \geq A\| p\|$; 
\item[$\bullet$] $C^2$-convex in the fiber i.e. for every $(q,p)\in T^*M$, the Hessian $\frac{\partial^2H}{\partial p^2}$ of $H$ in the fiber direction is positive definite as a quadratic form.\\
 \end{enumerate}
 We denote the Hamiltonian flow of $H$ by $(\varphi^H_t)$ and the Hamiltonian vector-field by $X_H$.  Note that the flow of a Tonelli Hamiltonian defined on $\A_n$ is not necessarily  a globally positive diffeomorphism. A geodesic flow is an example of a Tonelli flow.
 For example, the flat metric on $\T^n$ corresponds to the Tonelli Hamiltonian $H_0(q, p)=\frac{1}{2}\| p\|^2$ and its time-one flow is nothing but the map $f_0$ that we defined in example \ref{exa2}.\\

 At the end of the '80s, J.~Mather extended Aubry-Mather theory to the Tonelli Hamiltonians, introducing the concept of globally minimizing orbits and minimizing measures (see  \cite{Mat2} and \cite{Man1}).\\
  To explain that,  we  associate to any Tonelli Hamiltonian $H:T^*M\rightarrow \R$ its Lagrangian $L:TM\rightarrow \R$ that is dual to $H$ via the formula:
 $$\forall (q,v)\in TM, L(q,v)=\sup_{p\in T^*_qM} (p.v-H(q,p)).$$
 Then $L$ is as regular as $H$ is and is superlinear and $C^2$-convex in the fiber direction (see for example \cite{Fat1}). Moreover, we have:
 $$L(q,v)+H(q, p)=p.v\Longleftrightarrow v=\frac{\partial H}{\partial p}(q, p)\Longleftrightarrow p=\frac{\partial L}{\partial v}(q,v).$$
 If $\gamma:[\alpha, \beta]Ê\rightarrow M$ is an absolutely continuous arc, its {\em Lagrangian action} is then:
 $$A_L(\gamma)=\int_\alpha^\beta L(\gamma(t), \dot\gamma(t))dt.$$
 \subsection{Minimizing measures}
 \subsubsection{Case of the globally positive diffeomorphisms of the 2-dimensional annulus}\label{ssmintwistmap}
 We use the notations that were introduced in subsection \ref{ssgpd}. In the 2-dimensional case, J.~Mather and Aubry \& Le Daeron proved in \cite{ALD} and  \cite{Mat1} the existence of orbits $(q_i,p_i)_{i\in\Z}$ for $F$ that are {\em globally minimizing}. This means that for every $\ell\in\Z$ and every $k\geq 2$, $(q_{\ell+1}, \dots, q_{\ell+k-1})$ is minimizing the function $\Fc$ defined by:
 $$\Fc(q_{\ell+1}, \dots, q_{\ell+k-1})=\sum_{i=\ell+1}^{k} S(q_{i-1}, q_i).$$
 Then each of these orbits $(q_i, p_i)_{i\in\Z}$ is supported in the graph of a Lipschitz map defined on a closed subset of $\T$, and there exists a  bi-Lipschitz orientation preserving homeomorphism $h:\T\rightarrow\T$ such that $(q_i)_{i\in\Z}=(h^i(q_0))_{i\in\Z}$. Hence each of these orbits has a {\em rotation number}.\\
 Moreover, for each rotation number $\rho\in\R$, there exists a minimizing orbit that has this rotation number and there even exists a {\em minimizing measure}, i.e. an invariant measure whose support  is compact and filled by globally minimizing orbit, such that all the orbits contained in the support have the same rotation number $\rho$. These supports, which are Lipschitz graphs above a subset of $\T$,  are sometimes called {\em Aubry-Mather set}.\\
  In the following picture that concerns the so-called standard twist map, you can observe some invariant curves, some Cantor subsets and some periodic islands that must contain one periodic point.
   \begin{center}
\includegraphics[width=7cm]{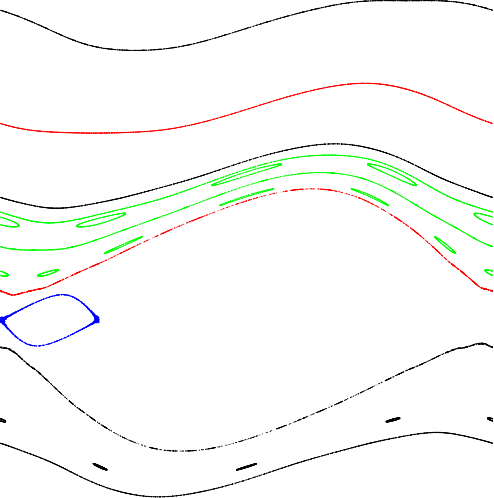}
\end{center} 

 Different kinds of Aubry-Mather sets can occur in this setting:
 \begin{enumerate}
 \item\label{AM1} some of them are invariant loops that are the graphs of some Lipschitz maps $\eta: \T\rightarrow \R$;
 \item \label{AM2} some other ones are just periodic orbits;
 \item\label{AM3} some of these Aubry-Mather sets are Cantor sets.
 \end{enumerate}
 In the case \ref{AM1}, it can happen that the dynamics restricted to the curve is bi-Lipschitz conjugate to a rotation;  in this case the Lyapunov exponents of the invariant measure supported in the curve are zero.  This is the case for the KAM curves. But P.~Le ~Calvez proved in \cite{LC1}  that in general (i.e. for a dense and $G_\delta$ subset of the set of the symplectic twist maps), there exists an open and dense subset $U$ of $\R$ such that any Aubry-Mather set that has its rotation number in $U$ is uniformly hyperbolic.

  \subsubsection{Case of the globally positive diffeomorphisms in higher dimension}
 
 For globally positive diffeomorphism in higher dimension, Garibaldi \& Thieullen prove the existence of globally minimizing orbits and measures in \cite{GarThi}. The results that they obtain are very similar to the ones that we recall in subsection \ref{ssminton} for Tonelli Hamiltonians.
 
 \begin{remk} There exists too an Aubry-Mather theory for time-one maps of time-dependent Tonelli Hamiltonians (see for example \cite{Ber1}). Even when the manifold $M$ is $\T^n$, the time-one map is not necessarily a  globally positive diffeomorphism of $\A_n$. Moreover, except for the 2-dimensional annulus (see \cite{Mos1}),  it is unknown if a globally positive diffeomorphism is always the time-one map of a time-dependent Tonelli Hamiltonian (see Theorem 41.1 in \cite{Gol1} for some partial results). In this survey, we won't speak about these time-one maps.
 
 \end{remk}
 
   \subsubsection{Case of the Tonelli Hamiltonians}\label{ssminton}
 It can be proved that is $q_b, q_e\in M$ are two points of $M$ and $\beta>\alpha$ two real numbers, if $\Gamma(q_b, q_e; \alpha, \beta)$ is the set of the $C^2$-arcs $\gamma: [\alpha, \beta]\rightarrow M$ that join $q_b$ to $q_e$ endowed with the $C^2$-topology, then  $\gamma$ is a critical point of the restriction of $A_L$ to $\Gamma(q_b, q_e; \alpha, \beta)$ if and only if $\gamma$ is the projection of an arc of orbit for $H$. This arc of orbit is then $(\gamma (t), \frac{\partial L}{\partial v}(\gamma(t), \dot\gamma(t)))_{t\in[\alpha, \beta]}$.
 
 In \cite{Mat2}, J.~Mather proves the existence of complete orbits $(\varphi_t^H(q,p))_{t\in\R}=(q(t), p(t))_{t\in\R}$ that are {\em globally minimizing}, i.e. such that every arc $(\pi\circ\varphi_t^H(x))_{t\in[\alpha, \beta]}=(q(t))_{t\in[\alpha, \beta]}$ is minimizing for the restriction of $A_L$ to $\Gamma(q(\alpha), q(\beta); \alpha, \beta)$.  He proves too the existence of {\em minimizing measures}, i.e. invariant measures whose support  is filled by globally minimizing orbit.
 
 When replacing $L$ by $L+\lambda$ where $\lambda$ is any closed 1-form on $M$, we obtain the same critical points for the Lagrangian action $A_{L+\lambda}$ as for the Lagrangian action $A_L$. But the minima for those two functions are not the same. Hence, adding different closed 1-form $\lambda$ to $L$ is a way to find other invariant measures supported in graphs, these measures being minimizing for $L+\lambda$.  The supports of these measures are the generalization of the Aubry-Mather sets. A rotation number can be associated to any minimizing measure (see \cite{Mat2}) and it can be proved that there exists a minimizing measure for any rotation number. But this doesn't give the existence of minimizing orbits of any rotation number (indeed the considered measures have not to be ergodic).

   \section{Lyapunov exponents for the minimizing measures and angle of the Oseledets splitting}\label{sLyap}
   Here we are interested in the Lyapunov exponents of the minimizing measures for globally positive diffeomorphisms or Tonelli Hamiltonian flows. In the case of symplectic twist maps, we noticed at the end of subsection \ref{ssmintwistmap} that these exponents may be non-zero or zero. 
   
  Let $(\Dc_t)$ ($t$ in $\Z$ or $\R$) be either the $\Z$-action generated by a globally positive diffeomorphism or the $\R$-action generated by a Tonelli Hamiltonian. Let $\mu$ be an ergodic  minimizing measure.  A general fact for ergodic measures and $C^1$-bounded dynamics is that the closer the stable and unstable bundles are (this means that there is an orthogonal basis of the stable bundle that is close to a orthonormal basis of the unstable bundle), the closer to zero the Lyapunov exponents are (see Proposition \ref{prop0}  in subsection \ref{sssLyap} for a more precise statement and \cite{Arn1} for a proof). 
    
    But in general, the converse assertion is false. We will see that it is true in the case of a twisting dynamics.  
    \begin{enumerate}
    \item[$\bullet$]  In subsection \ref{ssDirac} , we will prove these two statements for the Dirac measures  and even give a precise statement for the first assertion in the general case in   \ref{sssLyap}.
    \item[$\bullet$]  In subsection \ref{ssnbLyapunov} , we will explain that there is link between the number of non-zero Lyapunov exponents for a minimizing measure and the dimension of the intersection of the so-called Green bundles.
   \item[$\bullet$] In  subsection \ref{ssLyap}, we will explain the second assertion (and in fact a more precise statement): for a twisting dynamics with an hyperbolic minimizing measure, if   the stable and unstable bundles  are not close together (this means that all the unit vectors of the stable bundle are far from every unit vector of the unstable bundle), then all the positive Lyapunov exponents are large.
\end{enumerate}

   \subsection{Some simple remarks for Dirac masses} \label{ssDirac}
    Before looking at the Lyapunov exponents of any invariant measure, let us have a look to what happens for a Dirac mass in dimension 2. \\
    More precisely, let us assume that $x$ is a fixed point of a 2-dimensional diffeomorphism. We assume that $\sup\{ \| Df(x)\|, \| (Df(x))^{-1}\|\}\leq C$ where $C$ is some constant. Let $\lambda_1$, $\lambda_2$ be the two (complex) eigenvalues for $Df(x)$; then the Lyapunov exponents for the Dirac mass $\delta_x$ are $\log(|\lambda_1|)$ and $\log(|\lambda_2|)$. Let us assume that $\lambda_1$ and $\lambda_2$ are real and let us denote by $E_1$, $E_2$ the corresponding eigenspaces.
   \subsubsection{ What happens when the stable and unstable subspaces are close together} \label{sssLyap}
We have\\
   
 \noindent  {\em {\bf Simple principle: } if the eigenspaces $E_1$ and $E_2$ are close together, then the two eigenvalues $\lambda_1$ and $\lambda_2$ have to be close  together too.\\
   More precisely, if $e_i$ is a unit vector on $E_i$, we have:
   $$|\lambda_1-\lambda_2|\leq 2C\inf\{ \frac{\| e_2-e_1\|}{\| e_1+e_2\|}, \frac{\| e_1+e_2\|}{\| e_1-e_2\|}\}.$$
   }
  {\bf Proof of the simple principle.} We just compute
  $$Df(x).\frac{e_2-e_1}{\| e_2-e_1\|}=\frac{\lambda_2+\lambda_1}{2}.\frac{e_2-e_1}{\| e_2-e_1\|}+\frac{\lambda_2-\lambda_1}{2}\frac{e_2+e_1}{\| e_2-e_1\|}.$$
  As $e_2-e_1$ and $e_1+e_2$ are orthogonal, we deduce
  $$\frac{|\lambda_2-\lambda_1|}{2}\frac{\|e_2+e_1\|}{\| e_2-e_1\|}\leq \left\| Df(x).\frac{e_2-e_1}{\| e_2-e_1\|}\right\| \leq \| Df(x)\| \leq C.$$
  Changing $e_1$ into $-e_1$, we obtain the second inequality.\enddemo
  
   If $f$ is symplectic, we have $\lambda_2=\frac{1}{\lambda_1}$. In this case, if the eigenspaces are close together, the two eigenvalues have to be close to 1 and then the Lyapunov exponents are close to 0.
   
   This simple remark for fixed point can be generalized to any dimension and any invariant measure in the following way. For a proof, see \cite{Arn1}.
   
   \begin{nota}
We endow a compact manifold $N$ with a Riemannian metric. If $E$, $F$ are two linear subspaces of $T_xN$ that are $d$-dimensional with $d\geq 1$, the distance between $E$ and $F$ is:  
$${\rm dist} (E,F)=  \inf_{(e_i), (f_i)}\max\{ \| e_1-f_1\|, \dots , \| e_d-f_d\|\}  $$
where the infimum is taken over all the orthonormal basis $(e_i)$ of $E$, $(f_i)$ of $F$.
\end{nota}

   \begin{prop}\label{prop0}
Let $K$ be a compact subset of a manifold $N$,  let $C>0$   be a real number. Then, for any $f\in {\rm Diff}^1(M)$ so that $\max \{ \| Df_{|K}\|, \| Df^{-1}_{|K}\|\}\leq C$, if $f$ has an invariant ergodic measure $\mu$ with support in $K$ such that the Oseledets stable  and unstable bundles $E^s$ and $E^u$ of $\mu$ have the same dimension $d$  , if we denote by $\Lambda_u$   the sum of the positive   Lyapunov exponents and by $\Lambda_s$ the sum of the negative Lyapunov exponents, then: 
$$0< \Lambda_u-\Lambda_s\leq  d\log\left(1+(C^2+1)\int {\rm dist}(E^u, E^s)d\mu \right) $$
where ${\rm dist}$ is the distance.
\end{prop}

If for example $f$ is a symplectic diffeomorphism of $T^*M$, then $E^u$ and $E^s$ have same dimension (see for example \cite{BocVia1}). We deduce from the above proposition that if the stable and unstable Oseledets bundles are close together, then all the Lyapunov exponents are close to $0$. This result is not very surprising and not specific to the twisting dynamics. What is more surprising and specific to the twisting dynamics will come in the next section.
   \subsubsection{What happens when the stable and unstable subspaces are far from each other}
  For general symplectic dynamics,  we can have simultaneously two eigenvalues that are   close together and two eigenspaces that are not close together. See  for example the linear isomorphism of $\R^2$ with matrix in the usual basis:
  $$M=\begin{pmatrix}1+\varepsilon &0\\ 0& \frac{1}{1+\varepsilon}
  \end{pmatrix}
  $$
  with $\varepsilon>0$ small enough.
  
  As noticed by J.-C.~Yoccoz, this cannot happen for the minimizing Dirac masses of a twist map of the 2-dimensional annulus. It can be proved that at the fixed point $x$ corresponding to such a minimizing Dirac mass $\delta_x$, the eigenvalues of $Df(x)$ are real. We denote them by $\lambda_1$, $\lambda_2$ and by $E_1$, $E_2$ the two corresponding eigenspaces and by $M=\begin{pmatrix} a & b\\ c & d \end{pmatrix}$ the matrix of $Df(x)$ in coordinates $(q,p)$.  Then we have\\

\noindent{\em {\bf Simple result} If the torsion $b$ is bounded from below by a positive number, if $E_1$ and $E_2$ are far from each other, then $|\lambda_2-\lambda_1|$ cannot be to small.\\
More precisely, if $\theta$ is the angle between $E_1$ and $E_2$, then we have
$$|b|\leq\sup\{ 1,  \left({\rm cotan}(\theta)\right)^2\} |\lambda_2-\lambda_1|.$$
}
   {\bf Proof of the simple result.} The angle between $E_1$ and $E_2$ being $\theta$, there exists a matrix $R$ of rotation such that if $P:=R\begin{pmatrix} 1&{\rm cot an}\theta\\ 0&1\end{pmatrix}$, then $P\begin{pmatrix}1\\ 0\end{pmatrix}=e_1$ and  $P\begin{pmatrix}0\\ 1\end{pmatrix}\in\R.e_2$. As $R$ is a matrix of rotation, the modulus of all the coefficients of $P=\begin{pmatrix}\alpha&\beta\\ \delta&\gamma\end{pmatrix}$ is less that $\sup\{ 1, |{\rm cotan}\theta|\}$. Moreover, we have:
   $$M=\begin{pmatrix} \gamma & -\beta\\ -\delta& \alpha\end{pmatrix}.\begin{pmatrix} \lambda_1&0\\ 0&\lambda_2\end{pmatrix}.\begin{pmatrix} \alpha&\beta\\ \delta&\gamma\end{pmatrix}=\begin{pmatrix} *&\delta.\gamma (\lambda_1-\lambda_2)\\ *&*\end{pmatrix}.$$
   We deduce that $b=\delta\beta(\lambda_1-\lambda_2)$ and the wanted result.\enddemo
   
   \subsection{Number of non-zero Lyapunov exponents}\label{ssnbLyapunov}    Before giving an estimation of the non-zero Lyapunov exponents, we will try to find how many they are. As the dynamics is symplectic, we know that the number of negative Lyapunov exponents is equal to the number of positive Lyapunov exponents and then the number of zero Lyapunov exponents is even (see \cite{BocVia1} for a proof).
   \subsubsection{The two Green bundles}
   The Green bundles are two Lagrangian bundles that are defined along the minimizing orbits. In general, they are measurable but  not continuous. Let us recall that a subspace $H$ of the symplectic space $T_x(T^*M)$ is Lagrangian if its dimension is $n$ and if the restriction of the symplectic form to $H$ vanishes: $\omega_{|H\times H}=0$.
   
   The Green bundles were introduced in the '50s by L.~W.~Green to give a proof of the 2-dimensional version of Hopf conjecture: a Riemannian metric of $\T^n$ with no conjugate points is flat.  Then P.~Foulon extended the construction to the Finsler metrics in \cite{Fou1} and G.~Contreras \& R.~Iturriaga built them for any Tonelli Hamiltonian in \cite{ConItu1}. The construction for the twist maps of the annulus, and more generally for the twist maps of $\T^n\times \R^n$ is due to M.~Bialy \& R.~MacKay (see \cite{BiaMac}).
   
   We will  recall here their precise definition and we  will give their main properties. Before this, let us recall that there exists a way to compare different Lagrangian subspaces of $T_x(T^*M)$ that are transverse to the vertical $V(x)$. We choose some coordinates $(q, p)$ as explained at the beginning of section \ref{sectiontwist} and we denote the linearized coordinates by $(\delta q, \delta p)$ of $T_x(T^*M)$. If $H_1$, $H_2$ are two Lagrangian subspaces of $T_x(T^*M)$ that are transverse to the vertical $V(x)$,  we can write them in coordinates $(\delta q, \delta p)$ as the graph of some  symmetric matrices $S_1$, $S_2$. We say that $L_1$ is under $L_2$ and write $L_1\leq L_2$ when $S_2-S_1$ is a positive semi-definite matrix. We say that $L_1$ is strictly under $L_2$ and write $L_1<L_2$ if $L_1\leq L_2$ and $L_1$ and $L_2$ are transverse. This is equivalent to say that $S_2-S_1$ is positive definite. It can be proved that this definition doesn't depend on the chart that we choose. For an equivalent but more intrinsic definition, see \cite{Arn2}.

   Along every minimizing orbit of a globally positive diffeomorphism $F: \T^n\times \R^n\rightarrow \T^n\times \R^n$  or a Tonelli Hamiltonian flow $H: T^*M\rightarrow \R$ that we will denote by $(\Dc_t)$ (with $t$ in $\Z$ or $\R$), we can define two Lagrangian bundles $G_-$ and $G_+$.

\begin{defin}
If the orbit of $x$ is minimizing, then the familly $(D\Dc_t.V(\Dc_{-t}x))_{t>0}$  is a decreasing family of Lagrangian subspaces that converges to $G_+(x)$ and the familly $(D\Dc_{-t}.V(\Dc_{t}x))_{t>0}$  is an increasing family of Lagrangian subspaces that converges to $G_-(x)$.
\end{defin}

   We recall now some properties of the Green bundles.   \begin{enumerate}
   \item[$\bullet$] they are transverse to the vertical and $G_-\leq G_+$;
   \item[$\bullet$] $G_-$ and $G_+$ are invariant by the linearized dynamics, i.e. $D\Dc_t.G_\pm=G_\pm\circ \Dc_t$;
   \item[$\bullet$] for every compact $K$ such that the orbit of every point of $K$ is  minimizing, the two Green bundles restricted to $K$ are uniformly far from the vertical;
   \item[$\bullet$] (dynamical criterion) if the orbit of $x$ is minimizing and relatively compact in $T^*M$, if $\displaystyle{\liminf_{t\rightarrow+\infty} \| D(\pi\circ \Dc_t)(x)v\|\leq +\infty}$ then $v\in G_-(x)$.\\
    If $\displaystyle{\liminf_{t\rightarrow+\infty} \| D(\pi\circ \Dc_{-t})(x)v\|\leq +\infty}$ then $v\in G_+(x)$.
   \end{enumerate}
The bundles $G_-$ and $G_+$ are the {\em Green bundles}. The proof of the results that we mentioned before can be found in \cite{Arn2} for the Tonelli Hamiltonians and in \cite{Arn1} for the globally positive diffeomorphisms.\\
An easy consequence of the dynamical criterion and the fact that the Green bundles are Lagrangian is that when there is a splitting of $T_x(T^*M)$ into the sum of a stable, a center and a unstable bundles $T_x(T^*M)=E^s(x)\oplus E^c(x)\oplus E^u(x)$, for example an Oseledets splitting or a partially hyperbolic splitting, then we have
$$E^s\subset G_-\subset E^s\oplus E^c\quad{\rm and} \quad E^u\subset G_+\subset E^u\oplus E^c.$$
Let us give the argument of the proof. Because of the dynamical criterion, we have $E^s\subset G_-$. Because the dynamical system is symplectic, the symplectic orthogonal subspace to $E^s$ is  $(E^s)^\bot=E^s\oplus E^c$ (see e.g. \cite{BocVia1}). Because $G_-$ is Lagrangian, we have $G_-^\bot=G_-$. We obtain then $G_-^\bot=G_-\subset E^{s\bot}=E^s\oplus E^c$.\\
Let us note the following straightforward consequence: for a minimizing measure,    the whole information concerning the positive (resp. negative) Lyapunov exponents is contained in the restricted linearized dynamics $D\Dc_{t|G_+}$ (resp.  $D\Dc_{t|G_-}$). In particular, when the measure is weakly hyperbolic, we have almost everywhere $G_+=E^u$ and $G_-=E^s$.

\begin{nota}
Using a Riemannian metric on $M$, we define the horizontal subspace $\Hc$ as the kernel of the connection map. Then, for every Lagrangian subspace $\Gc$ of $T_x(T^*M)$, there exists a linear map $G: \Hc(x)\rightarrow V(x)$ whose graph is $\Gc$. That is the meaning of {\em graph} in what follows. When $M=\T^n$, we choose of course $\Hc=\R^n\times \{ 0\}$.\\
We denote by $s_+$ (resp. $s_-$) the linear map $  \Hc \rightarrow V $ with graph $G_+$ (resp. $G_-$). When we use symplectic coordinates, their matrices are symmetric. \\
Along a minimizing orbit in the case of a globally positive diffeomorphism, $G_k=Df^k (V\circ f^{-k}) $ (resp. $G_{-k}=Df^{-k} (V\circ f^k) $)  is the graph of $s_k$ (resp. $s_{-k}$). 
\end{nota}

\subsubsection{ÊLink between Êthe central dimension and the dimension of $G_-\cap G_+$Ê}
From $E^s\subset G_-\subset E^s\oplus E^c$ and $E^u\subset G_+\subset E^u\oplus E^c$, we deduce that $G_-\cap G_+\subset E^c$. Hence $G_-\cap G_+$ is an isotropic subspace (for $\omega$) of the symplectic space $E^c$. We deduce that $\dim (E^c)\geq 2\dim (G_-\cap G_+)$. When $E^s\oplus E^c\oplus E^u$ designates the Oseledets splitting of some minimizing measure $\mu$ , what is proved in \cite{Arn3} is that this inequality is an equality $\mu$ almost everywhere  for the Tonelli Hamiltonian flows and the same result is proved for the globally positive diffeomorphisms in \cite{Arn4}.
\begin{thm}\label{T1} Let $(\Dc_t)$ ($t$ in $\Z$ or $\R$) be either the $\Z$-action generated by a globally positive diffeomorphism or the $\R$-action generated by a Tonelli Hamiltonian. Let $\mu$ be a minimizing measure and let us denote by $p$ the $\mu$-almost everywhere dimension of $G_-\cap G_+$. Then $\mu$ has exactly $2p$ zero Lyapunov exponents, $n-p$ positive Lyapunov exponents and $n-p$ negative Lyapunov exponents.
\end{thm}

The idea is the following one. Firstly, let us notice that we have nothing to prove when $\dim(G_-\cap G_+)=n$ because we know that $\dim(E^c)\geq 2\dim (G_-\cap G_+)=2n$; in this case, $\dim(E^c)=2n$ and all the Lyapunov exponents are zero.

In the other case, we consider the following restricted-reduced linearized dynamics.  Let $\mu$ be an ergodic minimizing measure. Then the quantity $\dim(G_-\cap G_+)$ is $\mu$-almost everywhere constant. We denote this dimension by $p$ and we assume that $p<n$. 
\begin{nota}
We introduce the following linear spaces (see \cite{Arn3}): $E=G_-+ G_+$, $R=G_-\cap G_+$, $F$ is the reduced space $F=E/R$. As $E$ is coisotropic for $\omega
$ with $E^{\bot\omega}=R$, then $F$ is the symplectic reduction of $E$. As $E$ and $R$ are invariant by the linearized dynamics, then we can define a cocycle $M_t$ on $F$ as the reduced linearized dynamics. This cocycle is then symplectic for the reduced symplectic form $\Omega$.
\end{nota}

In \cite{Arn3} and \cite{Arn4}, we define for the cocycle $(M_t)$ a vertical subspace, some reduced Green bundles $g_-$ and $g_+$ that have properties similar to the ones of $G_\pm$, and we prove that $g_-$ and $g_+$ are transverse $\mu$-almost everywhere. As we will explain in next subsection, the transversality of the Green bundles implies the (weak) hyperbolicity of the measure. Here we have only the transversality of the reduced Green bundles, but this imply that the cocycle $(M_t)$ is (weakly) hyperbolic and then that the linearized dynamics has at least $2(n-p)$ non-zero Lyapunov exponents. This gives the conclusion.

\subsubsection{The transversality of the two Green bundles implies some hyperbolicity}
We will explain here why a minimizing measure $\mu$ is weakly hyperbolic when the Green bundles are transverse almost everywhere. We will deal with the discrete case (i.e.  globally positive diffeomorphisms of $\A_n$).  The diffeomorphism is denoted by $f$ and we assume that $\mu$-almost everywhere we have: $T_x\A_n=G_-(x)\oplus G_+(x)$.  We want to prove that $f$ has at least $n$ positive Lyapunov exponents; in this case,  because $f$ is symplectic,  $\mu$ has also $n$ negative Lyapunov exponents (see \cite{BocVia1}).

The idea is to use a bounded (but non continuous) symplectic change of linearized coordinates along the minimizing orbits where $T_x\A_n=G_-(x)\oplus G_+(x)$ such that $G_+$ becomes the horizontal and that preserves the vertical space.  Because $G_+$ is invariant by $Df$, the symplectic matrix of $Df^k$ is: $M^k(x)=\begin{pmatrix}
 a_k(x)&b_k(x)\\
0&d_k(x)\\
 \end{pmatrix}$.
 

 Because $G_k$ is transverse to the vertical, we have $\det b_k\not=0$. Because of the definition of $G_k$, we have then $d_k(x)=s_k(f^kx)b_k(x)$. As $(s_k(x))_{k\geq 1}$ is decreasing and tends to ${\bf 0}$ (because the horizontal is $G_+$), the symmetric matrix $s_k(f^kx)$ is positive definite. Moreover, because the matrix $M^k(x)$ is symplectic, we have:
   $$\left(M^k(x)\right)^{-1}=\begin{pmatrix} {}^td_k(x)&-{}^tb_k(x)\\
 0&{}^ta_k(x)\\
\end{pmatrix}$$
and by definition of $G_{-k}(x)$: 
${}^ta_k(x)=-s_{-k}(x){}^tb_k(x)$ and finally we have
$$M^k(x)=\begin{pmatrix}
-b_k(x)s_{-k}(x)& b_k(x)\\
0& s_k(f^kx)b_k(x)\\
\end{pmatrix}.
$$
The proof is then made of several lemmata. The first one is  a consequence of Egorov theorem and of the fact that $\mu$-almost everywhere on ${\rm supp}\mu$, $G_+$ and $G_-$ are transverse and then $-s_-$ is positive definite.
\begin{lemma}\label{LJ}
For every $\varepsilon >0$, there exists a measurable subset $J_\varepsilon$ of ${\rm supp}\mu$ such that:
\begin{enumerate}
\item[$\bullet$] $\mu (J_\varepsilon)\geq 1-\varepsilon$;
\item[$\bullet$] on $J_\varepsilon$, $(s_k)_{k\geq 1} $ and $(s_{-k})_{k\geq 1}$ converge uniformly ;
\item[$\bullet$] there exists a constant  $\alpha=\alpha(\varepsilon)>0$ such that: $\forall x\in J_\varepsilon,  -s_-(x)\geq \alpha {\bf 1}$.
\end{enumerate}
\end{lemma}
 We deduce:
 \begin{lemma}\label{LCVU}
 Let $J_\varepsilon$ be as in the previous lemma. On the set $\{ (k,x)\in\N\times J_\varepsilon, f^k(x)\in J_\varepsilon\}$, the sequence of conorms $(m(b_k(x))$ converges uniformly to $+\infty$, where $m(b_k)=\| b_k^{-1}\| ^{-1}$.
 \end{lemma}
 \demo Let $k, x$ be as in the lemma.\\
 The matrix  $M_k(x)=\begin{pmatrix}
 -b_k(x)s_{-k}(x)& b_k(x)\\
 0& s_k(f^kx)b_k(x)\\
 \end{pmatrix}
 $ being symplectic, we have: \\
 $-s_{-k}(x){}^tb_k(x)s_k(f^kx)b_k(x)={\bf 1}$ and thus  
 $-b_k(x)s_{-k}(x){}^tb_k(x)s_k(f^kx)={\bf 1}$ and:\\
  $b_k(x)s_{-k}(x){}^tb_k(x)=-\left(s_k(f^kx)\right)^{-1}$. \\
 We know that on $J_\varepsilon$, $(s_k)$ converges uniformly to zero. Hence,  for every $\delta>0$, there exists $N=N(\delta) $ such that: $k\geq N\Rightarrow \| s_k(f^kx)\|\leq \delta$. 
 
Moreover, as $G_{-1}\leq G_{-k}\leq G_{1}$ and $G_{-1}$ and $G_{1}$ continuously depend on $x$ in the compact subset ${\rm supp}\mu$ and because the linear change of coordinates that we use is bounded, there exists $\beta>0$ so that  $\| s_{\pm k}\|\leq \beta$ uniformly in $k$ on ${\rm supp}\mu$. Hence, if we choose $\delta'=\frac{\delta^2}{\beta}$, for every $k\geq N=N(\delta')$ and $x\in J_\varepsilon$ such that $f^kx\in J_\varepsilon$, we obtain: 
 $$\forall v\in\R^p,\beta \| {}^tb_k(x)v\|^2= {}^tv b_k(x)(\beta{\bf 1}){}^tb_k(x)v\geq - {}^tv b_k(x)s_{-k}(x){}^tb_k(x)v={}^tv\left(s_k(f^kx)\right)^{-1}v$$
 and we have: ${}^tv\left(s_k(f^kx)\right)^{-1}v\geq \frac{\beta}{\delta^2}\| v\|^2$ because $s_k(f^kx)$ is a positive definite matrix that is less than $\frac{\delta^2}{\beta}{\bf 1}$. We finally obtain: $\| {}^tb_k(x)v\|\geq \frac{1}{\delta}\| v\|$ and then the result that we wanted.
 \enddemo
 
 \begin{nota}
 We choose $\beta>0$ as in the previous proof, i.e. such that: $\forall k\in \Z\backslash\{ 0\}, \forall x\in{\rm supp}\mu, \| s_k(x)\|\leq \beta$.
 \end{nota}
 
  From now we fix  a small constant $\varepsilon>0$, associate  a set $J_\varepsilon$ with $\varepsilon$  via Lemma \ref{LJ} and a constant  $0<\alpha<\beta$; then  there exists $N\geq 0$ such that
 $$\forall x\in J_\varepsilon, \forall k\geq N, f^k(x)\in J_\varepsilon\Rightarrow m(b_k(x))\geq \frac{2}{\alpha}.$$

 \begin{lemma}
 Let $J_\varepsilon$ be as in Lemma \ref{LJ}. For $\mu$-almost point $x$ in $J_\varepsilon$, there exists a sequence of integers $(j_k)=(j_k(x))$ tending to $+\infty$ such that: 
 $$\forall k\in \N, m(b_{j_k}(x)s_{-j_k}(x))\geq \left( 2^\frac{1-\varepsilon}{2N}\right)^{j_k}.$$
 \end{lemma} 
 \demo
 As $\mu$ is ergodic for $f$, we deduce from Birkhoff ergodic theorem that for almost every point $x\in J_\varepsilon$, we have:
 $$\lim_{\ell\rightarrow +\infty}\frac{1}{\ell}\sharp \{ 0 \leq k\leq \ell-1; f^k(x)\in J_\varepsilon\}=\mu (J_\varepsilon)\geq 1-\varepsilon.$$
 We introduce the notation: $N(\ell)=\sharp \{ 0 \leq k\leq \ell-1; f^k(x)\in J_\varepsilon\}$.\\
 For such an $x$ and every $\ell\in\N$, we find a number $n(\ell)$ of integers:
 $$0=k_1\leq k_1+N\leq k_2\leq k_2+N\leq k_3\leq k_3+N\leq \dots \leq k_{n(\ell)}\leq \ell$$
 such that $f^{k_i}(x)\in J_\varepsilon$ and $n(\ell)\geq [\frac{N(\ell)}{N}]\geq \frac{N(\ell)}{N}-1$. In particular, we have: $\frac{n(\ell)}{\ell}\geq\frac{1}{N}(\frac{N(\ell)}{\ell}-\frac{N}{\ell})$, the right term converging to $\frac{\mu (J_\varepsilon)}{N}\geq \frac{1-\varepsilon}{N}$ when $\ell$ tends to $+\infty$. Hence, for $\ell$ large enough, we find: $n(\ell)\geq 1+  \ell \frac{1-\varepsilon}{2N}$.\\ 
As $f^{k_i}(x)\in J_\varepsilon$ and $k_{i+1}-k_i\geq N$, we have:  $m(b_{k_{i+1}-k_i}(f^{k_i}(x)))\geq \frac{2}{\alpha}$. Moreover, we have: $s_{-(k_{i+1}-k_i)}(f^{k_i}x)\leq s_-(f^{k_i}x)\leq -\alpha    {\bf 1}$ then $m(s_{-(k_{i+1}-k_i)}(f^{k_i}x))\geq \alpha$; hence: 
$$m(b_{k_{i+1}-k_i}(f^{k_i}x)s_{-(k_{i+1}-k_i)}(f^{k_i}x))\geq 2.$$
But the matrix $-b_{k_{n(\ell)}}(x)s_{-k_{n(\ell)}}(x)$ is the product of $n(\ell)-1$
 such matrices. Hence:
 $$m(b_{k_{n(\ell)}}(x)s_{-k_{n(\ell)}}(x))\geq 2^{n(\ell)-1}\geq 2^{\ell\frac{1-\varepsilon}{2N}}\geq \left( 2^\frac{1-\varepsilon}{2N}\right)^{k_{n(\ell)}}.
 $$

 \enddemo

This implies that all the Lyapunov exponents of the restriction of  $Df$ to $G_+$ are greater than $\log\left(  2^\frac{1-\varepsilon}{2N}\right)>0$.

\subsection{Lower bounds for the positive Lyapunov exponents} \label{ssLyap}

 \begin{nota}
 For a positive semi-definite symmetric matrix $S$ that is not the zero matrix, we   denote by $q_+(S)$ its smallest positive eigenvalue. 

 \end{nota}
 
 \begin{thm}\label{th2}
 Let $\mu$ be an ergodic  minimizing measure of a globally positive  diffeomorphism of $\A_n$ that has at least one non-zero Lyapunov exponent.We denote the smallest  positive  Lyapunov exponent of $\mu$ by $\lambda (\mu)$ and  an upper bound for $\| s_1-s_{-1}\|$ above ${\rm supp} \mu$ by   $C$. Then we have: 
$$\lambda (\mu)\geq\frac{1}{2}\int \log\left( 1+\frac{1}{C}q_+((s_+-\esse) (x))\right) d\mu (x).$$

 \end{thm}
The proof of this result is given in \cite{Arn4}.  
There is a similar result for Tonelli Hamiltonians: 
\begin{thm}\label{th3}
Let $\mu$ be an ergodic   minimizing  measure  for a Tonelli Hamiltonian $H: T^*M\rightarrow \R$  and with at least one non zero Lyapunov exponent; then its  smallest positive Lyapunov exponent $\lambda(\mu)$ satisfies:
$\lambda(\mu) \geq \frac{1}{2}\int m(\frac{\partial^2H}{\partial p^2}).q_+(s_+ -\esse) d\mu$.
\end{thm}
The proof of Theorem \ref{th2} is a little long and involves some technical changes of bases. We prefer to give the proof of Theorem \ref{th3}, that is simpler and shorter. The first point is the following lemma: 
 \begin{lemma}\label{lem9} Let $H: T^*M\rightarrow \R$ be a Tonelli Hamiltonian. Let $(x_t)$ be a   minimizing orbit and let $U$ and $S$ be two Lagrangian bundles along this orbit that are invariant by the linearized Hamilton flow and transverse to the vertical. Let $\delta x_U\in U$ be an infinitesimal orbit contained in the bundle $U$ and let us denote by $\delta x_S$ the unique vector of $S$ such that $\delta x_U-\delta x_S\in V$ (hence $\delta x_S$ is {\em not} an infinitesimal orbit). Then:
$$\frac{d{}}{dt}(\omega(x_t)(\delta x_S(t), \delta x_U(t)))={}^t(\delta x_U(t)-\delta x_S(t))  \frac{\partial^2H}{\partial p^2} (x_t)(  \delta x_U(t)-\delta x_S(t))\geq 0.$$
\end{lemma}
\demo As the result that we want to prove is local, we can assume that we are in the domain of a dual chart and express all the things in the corresponding dual linearized coordinates.\\
We consider an invariant Lagrangian  linear bundle $G$ that is transverse to  the vertical along the orbit of $x=(q,p)$. We denote the symmetric matrix whose graph is $G$  by $G$ again. An infinitesimal orbit contained in this bundle satisfies: $\delta p = G \delta q$. We deduce from the linearized Hamilton equations (if we are along the orbit $(q(t), p(t))=x(t)$, $\dot G$ designates $\frac{d{ }}{dt} (G(x(t)))$) that:
$$\delta \dot q = (\frac{\partial^2 H}{\partial q\partial p}   +\frac{\partial^2 H}{\partial p^2}   G)\delta q; \quad \delta \dot p = (\dot G +G\frac{\partial^2 H}{\partial q\partial p} +G\frac{\partial^2 H}{\partial p^2} G)\delta q=-(\frac{\partial^2 H}{\partial q^2}   +\frac{\partial^2 H}{\partial p\partial q} G)\delta q.$$
We deduce from these equations the classical Ricatti equation (it is  given for example in \cite{ConItu1} for Tonelli Hamiltonians, but the reader can find the initial and simpler Ricatti equation given by Green in the case of geodesic flows in \cite{Green}):
$$\dot G +G\frac{\partial^2 H}{\partial p^2}G+G\frac{\partial^2 H}{\partial q\partial p}+\frac{\partial^2 H}{\partial p\partial q}G+\frac{\partial^2 H}{\partial p^2}=0.$$
Let us assume now that the graphs of the symmetric matrices $\U$ and $\esse$ are invariant by the linearized flow  along the same orbit. We denote by $(\delta q_U, \U \delta q_U)$ an infinitesimal orbit that is contained in the graph of $\U$. Then we have: 
$$\frac{d{}}{dt}({}^t\delta q_U(\U-\esse) \delta q_U)=2{}^t\delta q_U(\U-\esse ) \delta\dot q_U +{}^t\delta q_U(\dot \U -\dot \esse )\delta q_U
$$
$$
=2{}^t\delta q_U(\U-\esse )(\frac{\partial^2 H}{\partial q\partial p} +\frac{\partial^2 H}{\partial p^2} \U )\delta  q_U +{}^t\delta q_U( \esse \frac{\partial^2 H}{\partial p^2}\esse - \U \frac{\partial^2 H}{\partial p^2} \U +\esse \frac{\partial^2 H}{\partial q\partial p}+\frac{\partial^2 H}{\partial p\partial q}\esse -\U\frac{\partial^2 H}{\partial q\partial p}-\frac{\partial^2 H}{\partial p\partial q}\U)\delta q_U
$$
$$
={}^t\delta q_U(\U \frac{\partial^2 H}{\partial q\partial p} -\esse \frac{\partial^2 H}{\partial q\partial  p}+\U \frac{\partial^2 H}{\partial p^2}s_+ -2\esse \frac{\partial^2 H}{\partial p^2} \U +\esse \frac{\partial^2 H}{\partial p^2}\esse  +\frac{\partial^2 H}{\partial p\partial q}\esse -\frac{\partial^2 H}{\partial p\partial q}\U )\delta q_U
$$
$$
={}^t\delta q_U (\U -\esse)\frac{\partial^2 H}{\partial p^2} (\U-\esse) \delta q_U \geq 0.
$$
To finish the proof, we just need to notice that in coordinates: $\omega(\delta x_S, \delta x_U)=$
$$\omega (\delta x_U, \delta x_U-\delta x_S)={}^t(\delta q_U, \U \delta q_U)\begin{pmatrix} 0& {\bf 1}\\
-{\bf 1} & 0\\
\end{pmatrix} \begin{pmatrix} 0\\ (\U-\esse)\delta q_U\\
\end{pmatrix}={}^t\delta q_U (\U-\esse)\delta q_U.
$$
\enddemo

Let   $\mu$ be an ergodic minimizing Borel probability measure  for a Tonelli Hamiltonian $H: T^*M\rightarrow \R$  and with at least one non zero Lyapunov exponent; its support $K$ is  compact and then,   there exists a constant $C>0$ such that $s_+$ and $s_-$ are bounded by $C$ above $K$. We choose a point $(q, p)$ that is generic for $\mu$ and $\delta x_+=(\delta q, s_+ \delta q)$ in the   Oseledets bundle corresponding to the smallest positive Lyapunov exponent $\lambda(\mu)$ of $\mu$ and we introduce $\delta x_-=(\delta q, s_- \delta q)$. 
Using the linearized Hamilton equations (see Lemma \ref{lem9}), because $\omega(x_t)(\delta x_-, \delta x_+)={}^t\delta q(s_+ -s_-)\delta q$, we obtain:
 $$\frac{d{}}{dt}(({}^t\delta q(s_+ -s_-)\delta q)=
{}^t\delta q(s_+ -s_-)\frac{\partial^2H}{\partial p^2}(q_t, p_t)(s_+-s_-)\delta q.$$
 Let us notice that $(s_+-s_- )^\frac{1}{2}\delta q$ is contained in  the orthogonal space to the kernel of $s_+-s_-$. Hence:
  $$\frac{d{}}{dt}(({}^t\delta q(s_+ -s_-)\delta q)\geq m(\frac{\partial^2H}{\partial p^2})q_+(s_+-s_- ) {}^t\delta q(s_+ -s_-)\delta q.$$
Moreover $\delta q\notin \ker (s_+-s_-)$ because $(\delta q,s_+ \delta q)$ corresponds to a positive Lyapunov exponent and then $(\delta q, s_+ \delta q)\notin G_-\cap G_+$. Then :  
  $$\begin{matrix}\hfill\frac{2 }{T}\log (\|\delta q(T)\|)+&\frac{\log 2C}{T}\geq \frac{1}{T}\log ({}^t\delta q (T)(s_+-s_- ) (q_T, p_T)\delta q (T))\geq\hfill \\
 & \frac{1}{T}\log({}^t\delta q (0)(s_+-s_-)(q,p)\delta q (0))+\frac{1}{T}\int_0^Tm(\frac{\partial^2H}{\partial p^2}(q_t, p_t))q_+((s_+-s_-)(q_t, p_t))dt.\end{matrix}$$
Using Birkhoff's ergodic theorem, we obtain: 
  $$\lim_{T\rightarrow +\infty}\frac{1 }{T}\log (\|\delta q(T)\|)=\lambda(\mu)\geq\frac{1}{2}\int m(\frac{\partial^2H}{\partial p^2})q_+(s_+-s_-)d\mu.$$

   \section{Shape of the support of the minimizing measures and Lyapunov exponents}
   
   \subsection{Some notations and definitions}
   
 For any subset $A\not= \emptyset$ of a manifold $M$ and any point $a\in A$, different kinds of subsets of $T_aM$ can be defined, that are cones and  also a generalizations of the notion of tangent space to a submanifold. We introduce them here when $M=\R^n$, but by using some charts, the definition can be extended to any manifold.
 
 \begin{defin}
 Let $A\subset \R^n$ a non-empty subset of $\R^n$ and let $a\in A$ be a point of $A$. Then
 \begin{enumerate}
 \item[$\bullet$] the {\em contingent cone } to $A$ at $a$ is defined as being the set of all the limit points of the sequences $t_k(a_k-a)$ where $(t_k)$ is a sequence of real numbers and $(a_k)$ is a sequence of elements of $A$ that converges to $a$. This cone is denoted by $C_aA$ and it is a subset of $T_a \R^n$;
 \item[$\bullet$] the {\em limit contingent cone} to $A$ at $a$ is the set of the limit points of sequences $v_k\in C_{a_k}A$ where $(a_k)$ is any sequence of points of $A$ that converges to $a$. It is denoted by $\widetilde C_aA$ and it is a subset of $T_a \R^n$;
 \item[$\bullet$] the   {\em paratangent cone} to $A$ at $a$ is the set of the limit points of the sequences  
$$ \lim_{ k\rightarrow \infty} t_k(x_k-y_k)$$
where $(x_k)$ and $(y_y)$ are sequences of elements of $A$ converging to $a$ and  $(t_k)$ is a sequence of elements of $\R$. It is denoted by $P_aA$ and it is a subset of $T_a \R^n$.
 \end{enumerate} 
 \end{defin}
  The following inclusions are always satisfied
  $$C_aA\subset \widetilde C_aA\subset P_aA.$$ 
  Let us give an example of a contingent and paratangent cone at a point where $A$ has an angle.
     \begin{center}
\includegraphics[width=6cm]{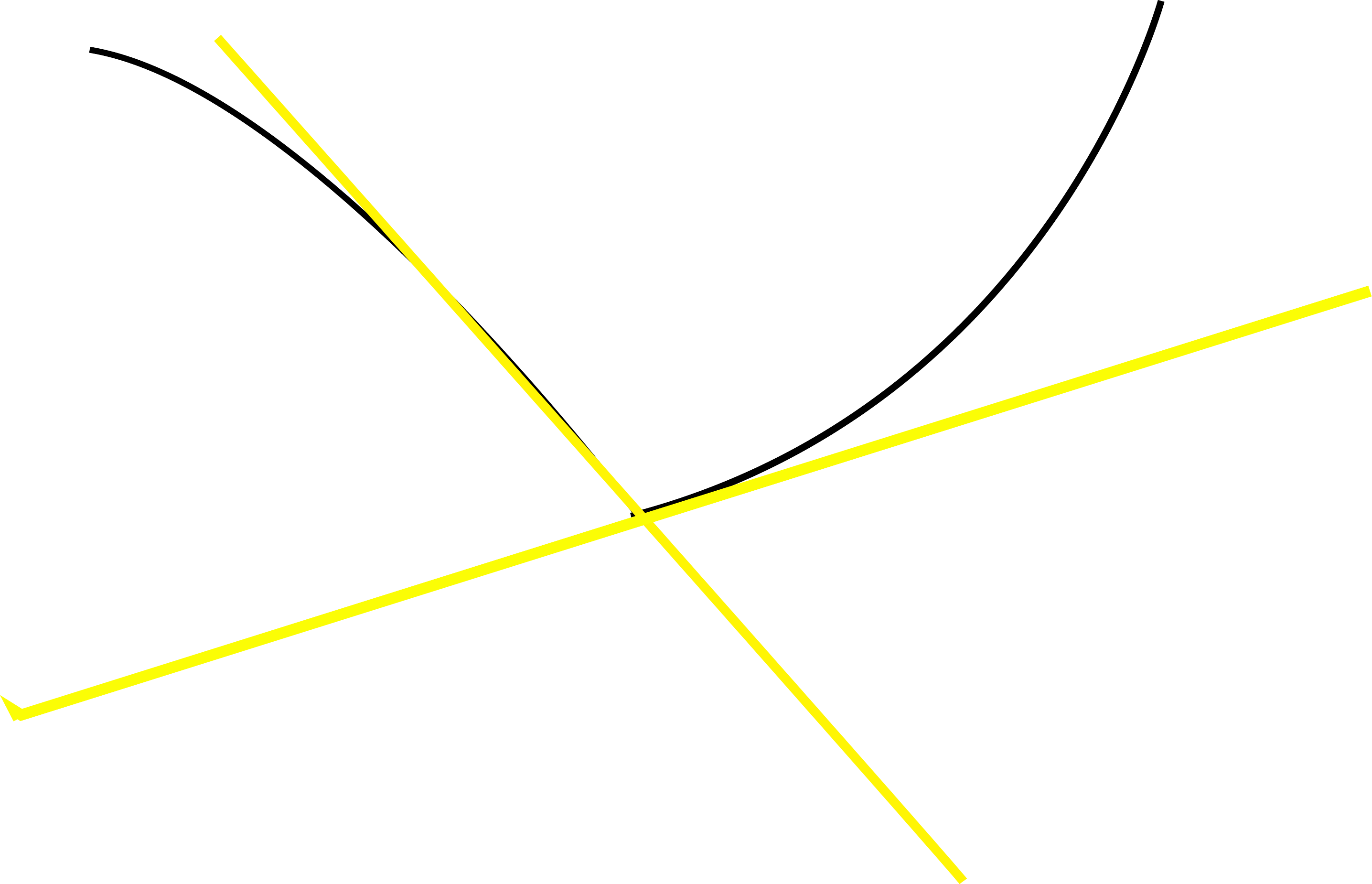}\includegraphics[width=6cm]{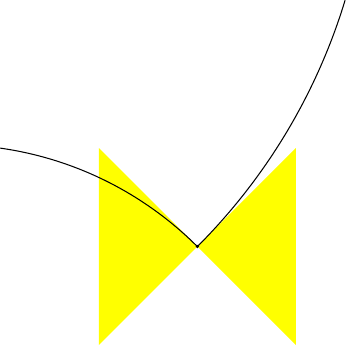}
\end{center} 
In the three last subsections of this survey, we will try to explain some relations between the Green bundles and these tangent cones. Unfortunately, in some cases, we need to use some modified Green bundles (see subsection \ref{ssmodified}). 
In general, the tangent cones are not Lagrangian subspaces (they are neither subspaces   nor isotropic). Because we need to compare them to Lagrangian subspaces, we give a definition:

\begin{defin}
Let $\Lc_-\leq \Lc_+$ be two Lagrangian subspaces of $T_x( T^*M) $ that are transverse to the vertical. If $v\in T_x (T^*M )$ is a vector, we say that $v$ is between $\Lc_-$ and $\Lc_+$ and write $\Lc_-\leq v\leq \Lc_+$ if there exists a third Lagrangian subspace in $T_x( T^*M) $ such that:
\begin{enumerate}
\item[$\bullet$] $v\in\Lc$;
\item[$\bullet$] $\Lc_-\leq \Lc\leq \Lc_+$.
\end{enumerate}
A subset $B$ of $T_x( T^*M) $ is between $\Lc_-$ and $\Lc_+$ if $\forall v\in  B, \Lc_-\leq v\leq \Lc_+$. Then we write $\Lc_-\leq B\leq \Lc_+$.
\end{defin} 
\begin{remk}
In the 2-dimensional case, $v$ is between $\Lc_-$ and $\Lc_+$ if and only if the slope of the line generated by $v$ is between the slopes of $\Lc_-$ and $\Lc_+$. In higher dimension, it is more complicated.
\end{remk}

\begin{defin}\begin{enumerate}
\item[$\bullet$] A subset $A$ of $\R^n\times \R^n$ is {\em $C^1$-isotropic} at some point $a\in A$ if $\widetilde C_aA$ is contained in some Lagrangian subspace;
\item[$\bullet$] a subset $A$ of $\R^n\times \R^n$ is {\em $C^1$-regular} at some point $a\in A$ if $P_aA$ is contained in some Lagrangian subspace.
\end{enumerate}
\end{defin}
Of course, the $C^1$-regularity of $A$ at a point $a$ implies the $C^1$-isotropy at the same point. But the converse implication is not true.\\
Observe that a $C^1$ Lagrangian submanifold is always $C^1$-regular.

   \subsection{Case I: 2-dimensional symplectic twist maps}
   The results that we explain now for the symplectic twist maps of the 2-dimensional annulus are proved in \cite{Arn0}.
   \begin{thm}\label{T4}
   Let $A$ be an Aubry-Mather set of a symplectic twist map of the 2-dimensional annulus $\A_1$. Then we have
   $$\forall a\in A, G_-(a)\leq P_aA\leq G_+(a).$$
   \end{thm}
   
    \begin{center}
\includegraphics[width=13cm]{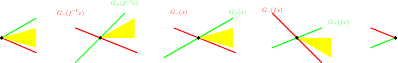}
\end{center} 
\begin{cor}\label{C1}
Let  $\mu$ be a minimizing ergodic measure of a symplectic twist map of the 2-dimensional annulus. If the Lyapunov exponents of $\mu$ are zero, then the support ${\rm supp}(\mu)$ of $\mu$ is $C^1$-regular $\mu$-almost everywhere.
\end{cor}
\begin{ques} Is there an example of such an invariant measure with zero Lyapunov exponents such that ${\rm supp}\mu$ is not $C^1$ at every point of ${\rm supp}\mu$?
\end{ques}
\begin{ques}\label{QUH} 
Is there an example of such an invariant measure with non-zero Lyapunov exponents such that ${\rm supp}\mu$ is not uniformly hyperbolic?
\end{ques}
Moreover, the following result is also true.
 \begin{prop}\label{P2}
Let  $\mu$ be a minimizing ergodic measure of a symplectic twist map of the 2-dimensional annulus that has an irrational rotation number. If the Lyapunov exponents of $\mu$ are non-zero, then the support ${\rm supp}(\mu)$ of $\mu$ is $C^1$-irregular $\mu$-almost everywhere.
 \end{prop}
\noindent We have even
  \begin{prop} 
Let  $\mu$ be a minimizing ergodic measure of a symplectic twist map of the 2-dimensional annulus that has an irrational rotation number. If the support ${\rm supp}(\mu)$ of $\mu$ is $C^1$-irregular   everywhere, then ${\rm supp}\mu$ is uniformly hyperbolic.
 \end{prop}
 Hence the size of the Lyapunov exponents can be read on the shape of ${\rm supp}(\mu)$. But how can we see in practice this irregularity? For example, if we want to ``draw'' (with a computer) our irregular (and hyperbolic) Aubry-Mather sets, we can use some sequences of minimizing periodic orbits. But if we look at the pictures of Aubry-Mather sets that exist, we see Cantor sets or curves, but we never see angles of the tangent spaces. That's why the following question was raised by X.~Buff~:\\
\begin{ques}(X.~Buff)  Is it possible (for example by using minimizing periodic orbits) to draw some Aubry-Mather sets with ``corners''? 
\end{ques}
   \subsection{Case II:  invariant Lagrangian graphs of Tonelli Hamiltonians}
   The proofs of the results we present in this section are given in \cite{Arn2}. We obtain a statement similar to Theorem \ref{T4} and Corollary \ref{C1} but no analogue to Proposition \ref{P2}.  Indeed, let us consider the following example: $(\psi_t)$ is a geodesic Anosov   flow defined on the cotangent bundle $T^*\Sc$ of a closed surface $\Sc$. Let $\Nc=T^*_1\Sc$ be its unit cotangent bundle, which is a 3-manifold invariant by $(\psi_t)$. Then a method due to Ma\~n\'e (see \cite{Man2}) allows us to define a Tonelli Hamiltonian $H$ on $T^*\Nc$ such that the restriction of its flow $(\varphi_t)$ to the zero section $\Nc$ is $(\psi_t)$: the Lagrangian $L$ associated with $H$  is defined by: $L(q,v)=\frac{1}{2}\| \dot\psi (q)-v\|^2$ where $\| .\|$ is any Riemannian metric on $\Nc$. In this case, the zero section is very regular (even $C^\infty$), but the Lyapunov exponents of every invariant measure with support  in $\Nc$ are non zero (except two, the one corresponding to the flow direction and the one corresponding to the energy direction). Hence, it may happen that some exponents are non zero and the support of the measure is very regular.
    \begin{thm}\label{T5}
   Let $\Gc$ be a Lipschitz Lagrangian   graph that is invariant by the flow of a Tonelli Hamiltonian $H:T^*M\rightarrow \R$. Then we have:
   $$\forall x\in \Gc, G_-(x)\leq P_x\Gc\leq G_+(x).$$
   \end{thm}
   The  following corollary is not proved in \cite{Arn2} but is an easy consequence of Theorems \ref{T5} and \ref{T1}.
   \begin{cor}\label{C2}
Let  $\mu$ be a minimizing ergodic measure for a Tonelli Hamiltonian of $T^*M$. If the Lyapunov exponents of $\mu$ are zero and if the support of $\mu$ is a graph above the whole manifold $M$, then the support ${\rm supp}(\mu)$ of $\mu$ is $C^1$-regular $\mu$-almost everywhere.
\end{cor}
\begin{ques} Is there an example where such a $\mu$ has zero Lyapunov and its support is not $C^1$   at least one point? \end{ques}
With the hypotheses of Corollary \ref{C2}, if we have further information about  the restricted dynamics to ${\rm supp}(\mu)$, we can improve the result in the following way.
\begin{prop}
 Let $\Gc$ be a Lipschitz Lagrangian  graph that is invariant  by the flow of a Tonelli Hamiltonian $H:T^*\T^n\rightarrow \R$. We assume that for some $T>0$, the restricted time-$T$ map $\varphi_{T|\Gc}$ is Lipschitz conjugated to some rotation of $\T^n$. Then $\Gc$ is the graph of a $C^1$ function.
 \end{prop}
 
 When $\mu$ is a minimizing measure with a support smaller than a Lagrangian graph, we don't obtain such a result (even if we have the feeling that it could be true).
 A fundamental tool  to  prove the previous results is the following proposition (that is proved in \cite{Arn2}).
 
 \begin{prop}\label{P4}
 Assume that the orbit of $x\in T^*M$ is globally minimizing for the Tonelli Hamiltonian $H: T^*M\rightarrow \R$ and that $\Lc$ defined on $\R$  is such that
 \begin{enumerate}
 \item[$\bullet$]  every $\Lc(t)$ is a Lagrangian subspace of $T_{\varphi_t^H(x)}(T^*M)$ that is transverse to the vertical subbundle;
 \item[$\bullet$] $\forall s, t\in\R, D\varphi_{t-s}^H\Lc(s)=\Lc(t)$.
 \end{enumerate}
 Then we have $\forall t\in\R, G_-(\varphi_t^H(x))\leq \Lc(t)\leq G_+(\varphi_t^H(x)).$
 \end{prop}
 Using Proposition \ref{P4} at any point where the invariant Lagrangian graph $\Gc$ is differentiable, we deduce a similar inequality for $\Lc$ being the tangent subspace at such a point. Then using a limit (and the notion of Clarke subdifferential), we deduce Theorem \ref{T5}.
 
 If we could obtain a result similar to Proposition \ref{P4} for vectors (instead of Lagrangian subspaces), we could deduce a similar statement for all minimizing measures. Hence we raise the question
 \begin{ques}\label{Q3}
Let $(\Dc_t)$ ($t$ in $\Z$ or $\R$) be either the $\Z$-action generated by a globally positive diffeomorphism or the $\R$-action generated by a Tonelli Hamiltonian.  Assume that the orbit of $x\in T^*M$ is globally minimizing and that the vector $v\in T_x(T^*M)$ is such that: $\forall t, D\varphi_t(v)\notin V(\Dc_tx)$. Is it true that:
 $$G_-(x)\leq v\leq G_+(x)?$$
 \end{ques}
 
 \begin{remk}
 Without a lot of change, all the results of this subsection  could be proved for any Lipschitz Lagrangian graph that is invariant by a globally positive diffeomorphism of $\T^n\times \R^n$.
 \end{remk}

  \subsection{Case III:  Tonelli Hamiltonians and globally positive diffeomorphisms}\label{ssmodified}
  The results contained in this subsection come from \cite{Arn3} and \cite{Arn4}. They use in a fundamental way a recent theory called the weak KAM theory  that was developped by  A.~Fathi in \cite{Fat1} in the case of the Tonelli Hamiltonians and by E.~Garibaldi \& P.~Thieullen in \cite{GarThi} in the case of the globally positive diffeomorphisms.
  
  Let us now introduce the modified Green bundles that we will use in this section. We use the constant $c_0=\frac{\sqrt{13}}{3}-\frac{5}{6}$. We identify $T_x(T^*M)$ to $\R^n\times \R^n$ in such a way that $\{0\}\times \R^n=V(x)$ is the vertical subspace and $\R^n\times\{ 0\}$ is the horizontal subspace $\Hc$.

\begin{defin}
We denote by $S_\pm(x):\R^n\rightarrow \R^n$ the linear operator such that $G_\pm(x)$ is the graph of $S_\pm(x)$: $G_\pm(x)=\{ (v, S_\pm(x)v); v\in\R^n\}$. Then the {\em modified Green bundles} $G_\pm$ are defined by:
$$\widetilde G_-(x)=\{ (v, (S_-(x)-c_0(S_+(x)-S_-(x)))v); v\in \R^n\}$${\rm and}$$ \widetilde G_+(x)=\{ (v, (S_+(x)+c_0(S_+(x)-S_-(x)))v); v\in \R^n\}.$$
\end{defin}
\begin{remk}
We have:
$$\widetilde G_-\leq G_-\leq G_+\leq \widetilde G_+.$$
Moreover, only the two following cases are possible
\begin{enumerate}
\item [$\bullet$] either $\widetilde G_-(x)$, $G_-(x)$, $G_+(x)$, $\widetilde G_+(x)$ are all distinct;
\item[$\bullet$] or $\widetilde G_-= G_-= G_+= \widetilde G_+$.
\end{enumerate}
\end{remk}

  \begin{thm}\label{T6} Let $\mu$ be a  minimizing measure for a Tonelli Hamiltonian of $T^*M$.  Then
$$\forall x\in {\rm supp} \mu, \widetilde G_-(x)\leq \widetilde C_x({\rm supp}\mu)\leq \widetilde G_+(x).$$
\end{thm}

Hence, the more irregular ${\rm supp}\mu$ is, i.e. the bigger the limit contingent cone is, the more distant  $\widetilde G_-$ and $\widetilde G_+$ (and thus $G_-$ and $G_+$ too)  are from each other and the larger the positive Lyapunov exponents are.

\begin{cor}\label{C3}
Let $H~: T^*M\rightarrow \R$ be a Tonelli Hamiltonian and let $\mu$ be an ergodic minimizing probability all of  whose   Lyapunov exponents are zero. Then, at $\mu$-almost every point of the support ${\rm supp}(\mu)$ of $\mu$, the set ${\rm supp}(\mu)$ is $C^1$-isotropic.
\end{cor}

There are two natural questions, that are related to question \ref{Q3} and that concern also the globally positive diffeomorphims.
\begin{ques}\label{Q4}
Can we replace $ \widetilde C_x({\rm supp}\mu)$ by $P_x({\rm supp}\mu)$ in Theorem \ref{T6} and Theorem \ref{Tcone}?
\end{ques}
\begin{ques}\label{Q5}
Can we replace $ \widetilde G_\pm(x)$ by $G_\pm(x)$ in Theorem \ref{T6} and Theorem \ref{Tcone}?
\end{ques}
 If the answer to question \ref{Q4} is positive, we can replace ``$C^1$-isotropic'' by ``$C^1$-regular'' in Corollary \ref{C3} and Corollary \ref{Cisotropic}.
 
 For the globally positive diffeomorphisms, we obtain a result only for the so-called strongly minimizing measures (the point is that for Tonelli Hamiltonians, miniminizing measures are also strongly minimizing).
 
 \begin{defin}
Let $F$ be a lift of a  {\em globally positive} diffeomorphism $f$ with generating function $S: \R^n\times \R^n\rightarrow \R$. A invariant Borel probability  $\nu$ on $\T^n\times \T^n$ is {\em strongly minimizing} if  $\nu$ is a minimizer in the following formula 
  $$\inf _\mu \int_{\R^n\times\R^n} S(x,y)d\tilde\mu(x,y);$$
  where the infimum is taken on the set of the Borel probability measures that are invariant by $f$ and $\tilde \mu$ is any lift of $\mu$ to a fundamental domain of $\R^n\times \R^n$ for the projection $(x,y)\mapsto (x, -\frac{\partial S}{\partial x}(x,y))$ onto $\T^n\times \R^n$.

 \end{defin} 
   E.~Garibaldi \& P.~Thieullen proved in \cite{GarThi} that such strongly minimizing measures exist. Moreover, they are minimizing.
   
   \begin{thm}\label{Tcone}
Let $\mu$ be a strongly minimizing measure of a globally positive  diffeomorphism of $\A_n$ et let ${\rm supp}\mu$ be its support. Then
$$\forall x\in {\rm supp} \mu, \widetilde G_-(x)\leq \widetilde C_x({\rm supp}\mu)\leq \widetilde G_+(x).$$

\end{thm}

\begin{cor}\label{Cisotropic}
Let $\mu$ be an ergodic strongly minimizing measure of a globally positive diffeomorphism of $\A_n$ all  of whose exponents are zero. Then ${\rm supp}\mu$ is $C^1$-isotropic almost everywhere.
\end{cor}

\end{document}